\documentclass{amsart}
\usepackage{amssymb}

\theoremstyle{plain}
\newtheorem{theorem}{Theorem}[section]
\newtheorem{lemma}[theorem]{Lemma}

\newtheorem{proposition}[theorem]{Proposition}
\newtheorem{corollary}[theorem]{Corollary}
\newtheorem{problem}[theorem]{Problem}
\newtheorem{observation}[theorem]{Observation}
\theoremstyle{definition}

\numberwithin{equation}{section}

\newcommand{\thmlabel}[1]{\label{thm:#1}} 
\newcommand{\lemlabel}[1]{\label{lem:#1}} 
\newcommand{\corlabel}[1]{\label{cor:#1}} 
\newcommand{\prplabel}[1]{\label{prp:#1}} 
\newcommand{\seclabel}[1]{\label{sec:#1}} 
\newcommand{\eqnlabel}[1]{\label{eqn:#1}} 

\newcommand{\thmref}[1]{\ref{thm:#1}}           
\newcommand{\lemref}[1]{\ref{lem:#1}}           
\newcommand{\corref}[1]{\ref{cor:#1}}           
\newcommand{\prpref}[1]{\ref{prp:#1}}           
\newcommand{\secref}[1]{\ref{sec:#1}}           
\newcommand{\eqnref}[1]{\eqref{eqn:#1}}     

\newcommand{\setof}[2]{\{ #1 \,|\, #2 \}}   
\newcommand{\sblp}[1]{\langle #1 \rangle}   

\newcommand{\Jordan}{\textsc{J}}
\newcommand{\jordan}[1]{\mathcal{J}(#1)}

\title{Admissible orders of Jordan loops}

\author[M. K. Kinyon]{Michael K.~Kinyon}
\author[K. Pula]{Kyle Pula}
\author[P. Vojt\v{e}chovsk\'y]{Petr Vojt\v{e}chovsk\'y}

\email[Kinyon]{mkinyon@math.du.edu}
\email[Pula]{jpula@math.du.edu}
\email[Vojt\v{e}chovsk\'y]{petr@math.du.edu}

\address{Department of Mathematics \\
University of Denver \\
2360 S Gaylord St \\
Denver, CO 80208, U.S.A.}

\begin{document}



\begin{abstract}
A commutative loop is Jordan if it satisfies the identity $x^2(yx) = (x^2y)x$.
Using an amalgam construction and its generalizations, we prove that a
nonassociative Jordan loop of order $n$ exists if and only if $n\ge 6$ and
$n\ne 9$. We also consider whether powers of elements in Jordan loops are
well-defined, and we construct an infinite family of finite simple
nonassociative Jordan loops.
\end{abstract}
\maketitle

\section{Introduction}

A magma $(Q,\cdot)$ is a \emph{quasigroup} if, for each $a,b\in Q$,
the equations $ax=b$, $ya=b$ have unique solutions $x,y\in Q$. A
\emph{loop} is a quasigroup with a neutral element. Standard
references on quasigroup and loop theory are \cite{Bruck,pflug}.

Most literature on loop theory focuses on varieties of loops satisfying some
near-associativity conditions with strong structural consequences. For
instance, \emph{Moufang loops} defined by the identity $x(y(xz)) = ((xy)x)z$
are \emph{diassociative} (every two elements generate a subgroup), and thus
\emph{power associative} (every element generates a subgroup).
In contrast, we will be dealing with Jordan loops whose properties are rather
weak indeed, and which consequently did not receive much attention yet.

A commutative loop is said to be \emph{Jordan} if it satisfies the \emph{Jordan
identity}
\begin{equation}\eqnlabel{Jordan}
    x^2(yx) = (x^2y)x   \tag{\Jordan}.
\end{equation}
This is the same identity used in the definition of
Jordan algebras \cite{McCrimmon}.

While studying Jordan loop rings, Goodaire and Keeping \cite{GoodaireKeeping}
asked for which integers $n$ there exists a nonassociative (that is, not
associative) Jordan loop of order $n$. We answer their question here:

\begin{theorem}\thmlabel{Main}
A nonasssociative Jordan loop of order $n$ exists if and only if $n\ge 6$ and
$n\ne 9$.
\end{theorem}

By the nature of the problem, this paper consists mostly of constructions, many
of which have a combinatorial flavor.

It is not difficult to find nonassociative Jordan loops of even order $n\ge 6$.
This part of the problem was already solved by Goodaire and Keeping in
\cite{GoodaireKeeping}. In \S\secref{Even}, we construct the same loops as in
\cite{GoodaireKeeping}, but using a visual argument. The core of our
construction is the classical correspondence between idempotent quasigroups and
loops of exponent $2$ (Proposition \prpref{Correspondence}).

The observation that an idempotent
quasigroup is a union of singleton subquasigroups suggests that the
aforementioned correspondence can be generalized. This is the idea
behind an amalgam construction of Foguel \cite{Foguel}:
idempotent quasigroups are generalized to quasigroups which are unions
of subquasigroups all of the same order. In \S\secref{Amalgam}, we use
this idea to build Jordan loops of all odd orders except
those of the form $2^m + 1$.

The amalgam construction can be generalized even further, this time
by allowing quasigroups that are unions of subquasigroups which are not
necessarily of the same order. In \S\secref{Fermat}, we use this
further generalization to construct Jordan loops of the troublesome
orders $2^m + 1$ for $m > 3$.

In \S\secref{Powers}, we ask if powers $x^k$ in Jordan loops are well-defined,
where by well-defined we mean that $x^k$ has the same value regardless of how
it is parenthesized. All orders $0 \leq k \leq 5$ are well-defined (Lemma
\lemref{Fifth}), and this is enough to show that the only Jordan loop of order
$5$ is the cyclic group (Corollary \corref{assoc5}). Up to relying on a
computer search for the case $n = 9$, this will complete the proof of Theorem
\thmref{Main}. A complete human proof of that case will appear elsewhere
\cite{Pula}.

We also show that for every composite integer $mn$ where $m > 1$ and
$n \geq 3$ is odd, there exists a $1$-generated Jordan loop with an element
$c$ such that  all powers $c^k$ are well-defined for $0\leq k < mn$, but for
which $c^{mn}$ is not well-defined (Theorem \thmref{powers}). The prime and
power of $2$ cases are more delicate and still open.

Finally, in \S\secref{Simple}, we construct an infinite family of nonassociative
simple Jordan loops of order $2^m - 1$.

We are pleased to acknowledge the assistance of the finite model builder
Mace4 developed by McCune \cite{mace4}. We found many of our constructions
by examining small Jordan loops built by Mace4. In addition, in this
paper, we rely on our Mace4 findings to complete the obstinate case $n = 9$ of
Theorem \thmref{Main}.

\section{The even order case}
\seclabel{Even}

The even order case is easy to handle, partly thanks
to the following observation.

\begin{lemma}
\lemlabel{CommExp2}
Every commutative loop of exponent $2$ is a Jordan loop.
\end{lemma}

\begin{proof}
Since each $x^2=1$, the Jordan identity \eqnref{Jordan} holds trivially.
\end{proof}

Our construction begins with the following classical correspondence,
which we need only in the commutative case.

\begin{proposition}
\prplabel{Correspondence}
There is a one-to-one correspondence between
commutative loops of exponent $2$ and order $n > 1$, and
commutative idempotent quasigroups of order $n-1$.
\end{proposition}

\begin{proof}
Given a commutative idempotent quasigroup on $X=\{x_2,\dots,x_n\}$,
introduce a new element $1$, let $1\cdot x = x\cdot 1 = x$ for every
$x\in X\cup \{1\}$, set $x_i\cdot x_i = 1$, and leave $x_i\cdot x_j$
intact for $i\ne j$. Conversely, given a nontrivial commutative loop
of exponent $2$ on $\{1,x_2,\dots,x_n\}$, remove $1$, set $x_i x_i=x_i$,
and leave $x_i\cdot x_j$ intact for $i\ne j$.
\end{proof}


Perhaps the most prominent application of the correspondence of Proposition
\prpref{Correspondence} is that between Steiner quasigroups (or Steiner triple
systems) and Steiner loops, i.e., commutative loops satisfying $x(xy)=y$
\cite{ColbournRosa}. We note in passing that it is not trivial to prove that a
Steiner loop of order $n$ exists if and only if $n$ is congruent to $2$ or $4$
modulo $6$ \cite{Kirkman}.

We now use Proposition \prpref{Correspondence} to construct the loops of
exponent $2$ given by Goodaire and Keeping in \cite{GoodaireKeeping}. We give
a visual argument, leaving to the reader the arithmetic details
(easily extracted from \cite{GoodaireKeeping}). First we show that
this construction can handle no more than the even order case.

\begin{lemma}
\lemlabel{2Odd} Let $Q$ be a finite commutative quasigroup. Then $|Q|$ is odd
if and only if the squaring map $x\mapsto x^2$ is a bijection.
\end{lemma}

\begin{proof}
Let $M$ be a multiplication table for $Q$ arranged so that columns and rows
are labeled in the same order. For each $a\in Q$, let $d(a)$ denote the
number of times $a$ appears on the main diagonal of $M$.
By commutativity, each $a$ occurs the same number of times above the main
diagonal of $M$ as it does below the main diagonal. Since there are
$|Q|$ occurrences of $a$ in $M$, it follows that $|Q| - d(a)$ is even.

Now if $x\mapsto x^2$ is not bijective, then by the finiteness of $Q$,
the map is not surjective. Thus for some $a \in Q$, $d(a) = 0$, and
so $|Q|$ is even.

Conversely, if $x\mapsto x^2$ is a bijection, then for each $a\in Q$,
$d(a) = 1$. Therefore $|Q|$ is odd.
\end{proof}

\begin{corollary}
\corlabel{2Odd}
Let $Q$ be a finite commutative loop. If there exists $a\in Q$, $a\neq 1$
such that $a^2 = 1$, then $|Q|$ is even. In particular, every
nontrivial, finite, commutative loop of exponent $2$ has even order.
\end{corollary}

\begin{corollary}
\corlabel{2Idem} Every finite commutative idempotent quasigroup has odd order.
\end{corollary}

\begin{lemma}
\lemlabel{IdempQuasi} For each odd $n > 0$, there exists a commutative
idempotent quasigroup of order $n$.
\end{lemma}

\begin{proof}
We will construct a multiplication table for such a quasigroup defined on the
set $S = \{1,\dots,n\}$, $n$ odd. Place $1,\dots,n$ on the diagonal, in this
order, that is, define $i\cdot i = i$ for all $i\in S$. Then fill out the rest
of the table so that all antidiagonals (with wrap around) are constant. This
works since $n$ is odd. The resulting quasigroup is clearly commutative and
idempotent.

Alternatively, if we let $S = \mathbb{Z}_n$ and define the quasigroup operation
$\ast$ by $a \ast b = (a + b)/2 \mod n$, we arrive at the same idempotent
quasigroup described above.
\end{proof}

Recall that a loop is \emph{left alternative} if it satisfies $x(xy) = (xx)y$.

\begin{proposition}
Let $n\ge 6$ be even. Then there is a commutative loop of exponent $2$
and order $n$ that is not left alternative, hence not associative.
\end{proposition}

\begin{proof}
Let $n\ge 6$ be even. Construct the commutative idempotent quasigroup $Q$ on
$\{2,\dots,n\}$ as in the proof of Lemma \lemref{IdempQuasi}. Let $L$ be the
commutative loop of exponent $2$ obtained from $Q$ by the correspondence of
Proposition \prpref{Correspondence}. We claim that $L$ is not left alternative.

Let $k = (n+2)/2$ be the average of $2$, $\dots$, $n$. Hence $k$ is on the main
antidiagonal in the canonical multiplication table of $Q$. Then $k\cdot (k+1) =
2$ and $k\cdot (k+2) = k+1$ (here we use the fact that $k+2\le n$, or,
equivalently, $n\ge 6$). Thus $(k\cdot k)\cdot (k+1) = 1\cdot (k+1) = k+1$,
while $k\cdot (k\cdot (k+1)) = k\cdot 2 \ne k+1$.
\end{proof}

It is well-known that every loop of order less than $5$ is associative. We have
therefore established the even case of Theorem \thmref{Main}.

\section{The amalgam construction}
\seclabel{Amalgam}

Let $(G,\circ)$ be a quasigroup, and let $S$ be a nonempty set.
Let $\triangledown = \setof{\triangledown_{g,h}}{g,h\in G}$ be
a family of binary operations on $S$ such that each
$Q_{g,h} = (S,\triangledown_{g,h})$ is a quasigroup, and
let $\mathcal{Q} = \setof{Q_{g,h}}{g,h\in G}$ denote the family
of all such quasigroups.
The \emph{quasigroup amalgam} $\mathcal{A}(G,\mathcal{Q})$
is the quasigroup defined on $S\times G$ by
\begin{equation}
\eqnlabel{QuasigroupAmalgam}
    (s,g)(t,h) = (s\triangledown_{g,h}t,g\circ h).
\end{equation}
Since it is always clear which multiplication $\triangledown_{g,h}$ is supposed
to be used in \eqnref{QuasigroupAmalgam}, we suppress the name of the binary
operation in $Q_{g,h}$ and write $st$ instead of $s\triangledown_{g,h}t$.

The term ``amalgam'' was suggested by Foguel \cite{Foguel}. In fact, this same
construction was first introduced by Bruck \cite[\S 10]{Bruck44}, who called it
an ``extension'' of $G$ by $S$. (However, this is not necessarily an extension
in the sense that the term is generally used in algebra.) An extensive list of
applications of the quasigroup amalgam can be found in \cite[pp. 36--43]{CPS}.

With the appropriate labeling of rows and columns, we can depict the
multiplication table of $\mathcal{A}(G,\mathcal{Q})$ as
\begin{displaymath}
    \begin{array}{c||c|c|c|c}
        \mathcal A&g&h&k&\dots\\
        \hline\hline
        &&&&\\
        g&(Q_{g,g},g\circ g)& (Q_{g,h},g\circ h) & (Q_{g,k},g\circ k) &\\
        &&&&\\
        \hline
        &&&&\\
        h&(Q_{h,g},h\circ g)&(Q_{h,h},h\circ h) & (Q_{h,k},h\circ k) &\\
        &&&&\\
        \hline
        &&&&\\
        k&(Q_{k,g},k\circ g)&(Q_{k,h},k\circ h) & (Q_{k,k},k\circ k) &\\
        &&&&\\
        \hline
        \vdots& & & &\ddots
    \end{array}
\end{displaymath}
where by $(Q_{g,h},g\circ h)$ we mean the $|S|\times |S|$ block
\begin{displaymath}
    \begin{array}{c|cccc}
    &(r,h)&(s,h)&(t,h)&\cdots\\
    \hline
    (r,g)&(rr,g\circ h)&(rs,g\circ h)&(rt,g\circ h)&\\
    (s,g)&(sr,g\circ h)&(ss,g\circ h)&(st,g\circ h)&\\
    (t,g)&(tr,g\circ h)&(ts,g\circ h)&(tt,g\circ h)&\\
    \vdots& & & &\ddots
    \end{array}
\end{displaymath}

We now attempt to turn the quasigroup amalgam into a loop by adjoining a new
neutral element $1$.

Firstly, we set $1x=x1=x$ for every $x\in (S\times G)\cup \{1\}$. Secondly, in
order to make sure that $1$ appears in every row and every column precisely
once, we select a bijection $c:G\to G$, for every $g\in G$ we select a loop
$L_g$ defined on $S\cup\{1\}$ with neutral element $1$, and we replace each
block $(Q_{g,c(g)},g\circ c(g))$ with the block $(L_g,g\circ c(g))$ from which
the row and column corresponding to $1$ has been removed. Finally, we identity
all elements of the form $(1,g)$ with $1$.

\begin{figure}
\begin{gather*}
    \begin{array}{|cc|cc|cc|}
    \hline
    (r,k)&(s,k)&\mathbf{(r,g)}&\mathbf{(s,g)}&(s,h)&(r,h)\\
    (s,k)&(r,k)&\mathbf{(s,g)}&\mathbf{(r,g)}&(r,h)&(s,h)\\
    \hline
    (s,g)&(r,g)&(s,h)&(r,h)&\mathbf{(s,k)}&\mathbf{(r,k)}\\
    (r,g)&(s,g)&(r,h)&(s,h)&\mathbf{(r,k)}&\mathbf{(s,k)}\\
    \hline
    \mathbf{(r,h)}&\mathbf{(s,h)}&(r,k)&(s,k)&(s,g)&(r,g)\\
    \mathbf{(s,h)}&\mathbf{(r,h)}&(s,k)&(r,k)&(r,g)&(s,g)\\
    \hline
    \end{array}
    \\
    \downarrow
    \\
    \begin{array}{|c|cc|cc|cc|}
    \hline
    1&(r,g)&(s,g)&(r,h)&(s,h)&(r,k)&(s,k)\\
    \hline
    (r,g)&(r,k)&(s,k)&\mathbf{(s,g)}&\mathbf{1}&(s,h)&(r,h)\\
    (s,g)&(s,k)&(r,k)&\mathbf{1}&\mathbf{(r,g)}&(r,h)&(s,h)\\
    \hline
    (r,h)&(s,g)&(r,g)&(s,h)&(r,h)&\mathbf{(s,k)}&\mathbf{1}\\
    (s,h)&(r,g)&(s,g)&(r,h)&(s,h)&\mathbf{1}&\mathbf{(r,k)}\\
    \hline
    (r,k)&\mathbf{(s,h)}&\mathbf{1}&(r,k)&(s,k)&(s,g)&(r,g)\\
    (s,k)&\mathbf{1}&\mathbf{(r,h)}&(s,k)&(r,k)&(r,g)&(s,g)\\
    \hline
    \end{array}
\end{gather*}
\caption{Attempting to create a loop out of a quasigroup amalgam by adjoining
$1$}\label{Fg:Process}
\end{figure}

The process is illustrated in Figure \ref{Fg:Process}, where the underlying
quasigroup $(G,\circ)$, the quasigroups $Q_{g,h}$ and the loops $L_g$ are
chosen as follows
\begin{displaymath}
G = \begin{array}{c|ccc}
    &g&h&k\\ \hline g&k&g&h\\ h&g&h&k\\ k&h&k&g
    \end{array},\quad
Q_{x,y} = \begin{array}{c|cc}
    &r&s\\ \hline r&r&s\\ s&s&r
    \end{array}\quad\text{or}\quad
    \begin{array}{c|cc}
    &r&s\\ \hline r&s&r\\ s&r&s
    \end{array},\quad
L_x = \begin{array}{c|ccc}
    &1&r&s\\ \hline 1&1&r&s\\ r&r&s&1\\ s&s&1&r
    \end{array},
\end{displaymath}
and where the bijection $c:G\to G$ is $c(g)= h$, $c(h)=k$, $c(k)=g$.

Although it might appear that the resulting quasigroup should be a loop, Figure
\ref{Fg:Process} shows otherwise. Here are the conditions that make the
construction work:

\begin{proposition}\prplabel{LoopAmalgam}
Let $G$, $L_g$, $Q_{g,h}$ and $c:G\to G$ be as above. Then the resulting
magma is a loop if and only if $c$ is the identity on $G$ and $G$ is an
idempotent quasigroup.
\end{proposition}

\begin{proof}
We show the necessity of the conditions and leave the sufficiency to the
reader. Assume that the resulting magma is a loop. Let $g\in G$ and $s\in S$.
Let $k\in G$ be such that $g\circ k = g$. If $c(g)\ne k$, then the block
$(Q_{g,k},g\circ k)$ has not been replaced, and thus it contains $(s,g\circ k)
=(s,g)$ in the row labeled by $(s,g)$, a contradiction with $(s,g)1=(s,g)$.
Hence $g\circ c(g)=g$. Suppose that $g\ne c(g)$. Then there is $h\ne g$ such
that $h\circ c(g)=c(g)$. Since the block $(Q_{h,c(g)},h\circ c(g))$ has not
been replaced, it contains $(s,h\circ c(g)) = (s,c(g))$ in the column labeled
by $(s,c(g))$, a contradiction with $1(s,c(g)) = (s,c(g))$. Thus $g=c(g)$.
\end{proof}

We therefore define the loop amalgam as follows: Let $(G,\circ)$ be an
idempotent quasigroup, $S$ a nonempty set, and $1$ an element not contained in
$S$. Let $\mathcal{Q} = \setof{Q_{g,h}}{g,h\in G, g\neq h}$ be a family of
quasigroups with underlying set $S$. Let $\mathcal{L} = \setof{L_g}{g\in G}$ be
a family of loops with underlying set $S\cup\{1\}$ and with neutral element
$1$. Then the (\emph{loop}) \emph{amalgam} $\mathcal{A} =
\mathcal{A}(G,\mathcal{L},\mathcal{Q})$ is a magma defined on the set $X =
(S\times G)\cup \{1\}$ by
\begin{equation}\eqnlabel{Amalgam}
    \begin{array}{c}
        1x = x1 = x\text{ for every $x\in X$},\\
        (s,g)(t,h) = (st,g\circ h)\text{ for every $s$, $t\in S$, $g$, $h\in G$},
    \end{array}
\end{equation}
where we identify the elements $(1,g)$ with $1$, and where the multiplication
$st$ takes place in $L_g$ when $g=h$, and in $Q_{g,h}$ when $g\ne h$.

The following result is then easily observed:

\begin{lemma}\
\lemlabel{CommSize}
The amalgam $\mathcal{A}(G,\mathcal{L},\mathcal{Q})$ is commutative if and
only if $(G,\circ)$ is commutative, each $L_g$ is commutative, and $Q_{g,h}$
is the opposite quasigroup of $Q_{h,g}$ for every $g\ne h\in G$. In addition,
$|\mathcal A| = |G||S| + 1 = |G|(|L_g|-1) + 1$.
\end{lemma}

Note that $\mathcal{Q} = \emptyset$ if $G = \{x\}$ is a singleton quasigroup.
Any loop $L$ can then be obtained as an amalgam by letting  $G=\{x\}$ and
$L_x = L$, a situation that we call \emph{trivial}.

A restricted version of the amalgam construction appears in \cite{Foguel}.
Namely, Foguel lets $(G,\circ)$ be an idempotent quasigroup, sets $L_g = L$
for every $g\in G$, and $Q_{g,h} = Q$ for every $g\ne h\in G$. We denote his
construction by $\mathcal{A}(G,L,Q)$.

\begin{proposition}
\prplabel{GuaranteedJordan}
Let $(G,\circ)$ be an idempotent quasigroup, $(L,\bullet)$
a loop, and $(Q,*)$ a quasigroup defined on $L\setminus \{1\}$. Then
$\mathcal{A} = \mathcal{A}(G,L,Q)$ is a Jordan loop if and only if
\begin{enumerate}
    \item $L$ is a Jordan loop, and
    \item $G$ and $Q$ are commutative, and
    \item for every $s$, $t\in Q$ either $s\bullet s=1$ or $(s\bullet
    s)*(t*s) = ((s\bullet s)*t)*s$.
\end{enumerate}
\end{proposition}

\begin{proof}
We know from Lemma \lemref{CommSize} that $\mathcal A$ is commutative if and
only if $G$, $L$ and $Q$ are commutative.

Let $x=(s,g)\ne 1\ne (t,h) = y\in\mathcal A$. If $g=h$ then \eqnref{Jordan}
holds for $x$, $y$ if and only if $(s\bullet s)\bullet(t\bullet s) = ((s\bullet
s)\bullet t)\bullet s$.

Assume that $g\ne h$. Then
\begin{align*}
    &x^2(yx) = (s\bullet s,g)(t*s,h\circ g),\\
    &(x^2y)x = (s\bullet s,g)(t,h)\cdot(s,g).
\end{align*}
If $s\bullet s=1$, the two expressions are equal. If $s\bullet s\ne 1$, then
\begin{align*}
    &x^2(yx) = ((s\bullet s)*(t*s),g\circ(h\circ g)),\\
    &(x^2y)x =((s\bullet s)*t,g\circ h)(s,g) = (((s\bullet s)*t)*s,(g\circ
    h)\circ g).
\end{align*}
Thus, when $G$ is commutative, the two expressions are equal if and only if
$(s\bullet s)*(t*s) = ((s\bullet s)*t)*s$.
\end{proof}

\begin{proposition}\prplabel{NearlyAllOdd}
Let $n>5$ be an odd integer not of the form $2^m + 1$. Then
there exists a Jordan loop of order $n$ which is not left
alternative, and hence, is not associative.
\end{proposition}

\begin{proof}
Since $n-1\ne 2^m$, we have $n-1 = 2^{\ell} k$ where $k\ge 3$ is odd and $\ell\geq 0$.
Let $(G,\circ)$ be a commutative idempotent quasigroup of order
$k$, which exists by Lemma \lemref{IdempQuasi}. Let $(L,\bullet)$ be a
commutative group of order $2^{\ell} + 1$, and $(Q,*)$ a commutative group
of order $2^{\ell}$ defined on $L\setminus\{1\}$.

Set $\mathcal{A} = \mathcal{A}(G,L,Q)$. Then $|\mathcal{A}| = n$ by
Lemma \lemref{CommSize}, and $\mathcal{A}$ is a Jordan loop by Proposition
\prpref{GuaranteedJordan}.

Choose distinct elements $g,h\in G$ and an element $s\in L$ with $s\neq 1$.
Then $g = g\circ g \ne g\circ h$, and since $L$ is a group of odd order,
$s\bullet s\ne 1$. Thus
\begin{align*}
    &(s,g)(s,g)\cdot (s,h) = (s\bullet s,g)(s,h) = ((s\bullet s)*s,g\circ h),\\
    &(s,g)\cdot (s,g)(s,h) = (s,g)(s*s,g\circ h) = (s*(s*s),g\circ(g\circ h)).
\end{align*}
As $g\ne g\circ h$, we have $g\circ h = (g\circ h)\circ(g\circ h)\ne
g\circ(g\circ h)$, and it follows that $\mathcal{A}$ is not left alternative.
\end{proof}

\section{Order $2^m+1$}\seclabel{Fermat}

Assume momentarily that $\mathcal{A}(G,\mathcal{L},\mathcal{Q})$ is a commutative
loop of order $2^m + 1$ and $|G| > 1$. Then, by Lemma \lemref{CommSize}, $G$ is a
commutative idempotent quasigroup of even order, which contradicts Lemma
\lemref{IdempQuasi}. It is therefore impossible to construct a commutative loop
of order $2^m + 1$ as a nontrivial amalgam.

When we arrived at this impasse, we seriously entertained the possibility that
there are infinitely many orders $n$ for which no nonassociative Jordan loops
exist. We soon noticed that this alternative has interesting number-theoretical
consequences. We were unable to resist the temptation to present them here,
despite the fact that Theorem \thmref{Main} renders them irrelevant.

\begin{observation}
\label{Ob:Fermat}
Assume that there are infinitely many orders $n$ for which no
nonassociative Jordan loop exists. Then there are infinitely many Fermat
primes, i.e., primes of the form $2^m+1$.
\end{observation}

\begin{proof}
The results of \S\secref{Even} and Proposition \prpref{NearlyAllOdd} imply
that under our assumption there are infinitely many orders $n=2^m+1$ for which
no nonassociative Jordan loop exists. Let $p_1^{a_1}\cdots p_k^{a_k}$ be a
prime factorization of such an $n=2^m+1$.

If some $p_i$ is not Fermat, then $p_i>5$, there is a nonassociative Jordan
loop $Q$ of order $p_i$ by Proposition \prpref{NearlyAllOdd}, and hence a
nonassociative Jordan loop of order $n$ is obtained as a direct product of $Q$
with a commutative group of order $n/p_i$, a contradiction. Thus every $p_i$ is
a Fermat prime.

If there is an $i$ such that $p_i>3$ and $a_i>1$ then $p_i^2 = (2^r+1)^2$ is
not of the form $2^s+1$, as can be routinely observed, so there is a
nonassociative Jordan loop $Q$ of order $p_i^2$ by Proposition
\prpref{NearlyAllOdd}, again leading to a contradiction. Hence for every $i$ we
have $p_i\le 3$ or $a_i=1$.

If there is an $i$ such that $p_i=3$ and $a_i>2$ then $p_i^3=27$ divides $n$
but is not of the form $2^s+1$, a contradiction once again.

Altogether, we see that every $p_i$ is a Fermat prime, if $p_i>3$ then $a_i=1$,
and if $p_i=3$ then $a_i\le 2$. Hence there must be infinitely many Fermat
primes.
\end{proof}

Since it is generally (but not universally) believed that (i) there are
finitely many Fermat primes, and (ii) it is hard to establish (i), we concluded
that the premise of Observation \ref{Ob:Fermat} is likely false. And, indeed, a
modification of the amalgam construction allowed us to settle the case $2^m + 1$,
too.

In \S\secref{Amalgam}, we showed how to turn the quasigroup amalgam
$\mathcal{A} = \mathcal{A}(G,\mathcal{Q})$ into a loop by adjoining a neutral element
$1$ and replacing the diagonal blocks $(Q_{g,g},g)$ by $(L_g,g)$. The
underlying principle that makes this work is the fact that $\mathcal A$ is a
quasigroup that is a disjoint union of (some of its) subquasigroups, namely
$\mathcal A = \bigcup_{g\in G} (Q_{g,g},g)$.

This leads us to a more general construction: Let $Q = \bigcup_{i\in I}Q_i$ be
a quasigroup, where $Q_i\le Q$ and $Q_i\cap Q_j=\emptyset$ for every $i\ne j\in
I$. Let $1$ be an element not contained in $Q$. For every $i\in I$, let $L_i$
be a loop defined on $Q_i\cup\{1\}$. Upon replacing the blocks $Q_i\times Q_i$
with $(L_i\setminus\{1\})\times(L_i\setminus\{1\})$ in the multiplication table
of $Q$ and upon adjoining $1$ as the neutral element, we obtain a loop.

Instead of introducing more notation, we say, somewhat informally, that the
resulting loop is obtained from $Q$ by \emph{replacing the subquasigroups} $Q_i$
\emph{with subloops} $L_i$.

Let us a call a quasigroup \emph{Jordan} if it is commutative and satisfies
\eqnref{Jordan}.

\begin{lemma}\lemlabel{JordanQuasigroup}
Let $G$ be a commutative group, $G_i\le G$ for $i\in I$, $G_i\cap G_j=\{1\}$
for $i\ne j\in I$, and $G=\bigcup_{i\in I}G_i$. Let $Q=G\setminus\{1\}$,
$Q_i=G_i\setminus\{1\}$, and assume that $(Q_i,*)$ is a Jordan quasigroup for
every $i\in I$. Then $(Q,*)$ with multiplication defined by
\begin{displaymath}
    x*y = \left\{\begin{array}{ll}
        x*y,&\text{ if $x$, $y\in Q_i$},\\
        xy,&\text{ if $x\in Q_i$, $y\in Q_j$, $i\ne j$}
        \end{array}\right.
\end{displaymath}
is a Jordan quasigroup, $Q_i\le Q$, $Q_i\cap Q_j=\emptyset$ for $i\ne j$, and
$Q=\bigcup_{i\in I}Q_i$.
\end{lemma}
\begin{proof}
We leave it to the reader to verify that $(Q,*)$ is a commutative quasigroup.
Let us check \eqnref{Jordan}. Let $x$, $y\in Q$. If there is $i\in I$ such
that $x\in Q_i$ and $y\in Q_i$, then the verification of $(x*x)*(y*x) =
((x*x)*y)*x$ takes place entirely within the Jordan quasigroup $(Q_i,*)$.
Assume that $x\in Q_i$, $y\in Q_j$, and $i\ne j$. Since $Q_iQ_j\cap
Q_i=\emptyset$ and $Q_i*Q_i\subseteq Q_i$, we have
\begin{align*}
    &(x*x)*(y*x) = (x*x)*(yx) = (x*x)(yx),\\
    &((x*x)*y)*x = ((x*x)y)*x = ((x*x)y)x.
\end{align*}
We are through because $G$ is associative.
\end{proof}

\begin{proposition}\prplabel{2mPlus1}
For every $m>3$ there is a nonassociative Jordan loop of order $2^m+1$.
\end{proposition}
\begin{proof}
Let $G = \langle
\alpha,\beta;\;\alpha^3=\beta^3=1,\,\alpha\beta=\beta\alpha\rangle$ be the
direct product of the cyclic group of order $3$ with itself. Let $G_1 = \langle
\alpha\rangle$, $G_2 = \langle \beta\rangle$, $G_3 = \langle
\alpha\beta\rangle$, $G_4 = \langle \alpha\beta^2\rangle$. Then $G_i\cap
G_j=\{1\}$ for $i\ne j$, and $G=\bigcup_{1\le i\le 4}G_i$.

For $Q_i=G_i\setminus\{1\}$, let $(Q_i,*)$ be the cyclic group of order $2$.
Let $(Q,*)$ be the Jordan quasigroup (of order $8$) obtained by Lemma
\lemref{JordanQuasigroup}.

Let $C$ be the cyclic group of order $2^{m-3}$, and set $(\overline Q,*) =
C\times Q$. Then $\overline{Q}$ is a Jordan quasigroup of order $2^m$ that is a
disjoint union of its subquasigroups $\overline{Q_i} = C\times Q_i$, each of
order $2^{m-2}$.

Let $L$ be a commutative group of order $2^{m-2}+1$, and let $J$ be the loop
obtained from $\overline{Q}$ by replacing the subquasigroups $\overline{Q_i}$
with copies of $L$. Clearly, $|J|=2^m+1$, $J$ is commutative (because $Q$ is),
and we claim that $J$ is a nonassociative Jordan loop.

We first check that $J$ satisfies the Jordan identity $x^2(yx)=(x^2y)x$. If
$x=1$ or $y=1$, there is nothing to prove. If $x$, $y$ belong to the same
subquasigroup $\overline{Q_i}$ of $\overline{Q}$, then the verification of
$x^2(yx)=(x^2y)x$ takes place entirely within a copy of the group $L$. Assume
that $x\in\overline{Q_i}$, $y\in\overline{Q_j}$, $i\ne j$. Since $L$ contains
no elements of order $2$, we have $x^2\in \overline{Q_i}$. Then $x^2(yx) =
x^2(y*x) = x^2*(y*x)$, and $(x^2y)x = (x^2*y)x = (x^2*y)*x$. As
$(\overline{Q},*)$ is Jordan, \eqnref{Jordan} holds for $J$, too.

Since $L$ is a subloop of $J$, $|L|=2^{m-2}+1$, $|J|=2^m+1$, and $2^{m-2}+1$
does not divide $2^m+1$ for $m>3$, it follows that $J$ cannot be a group,
by the Lagrange theorem.
\end{proof}

The restriction that $m > 3$, which is not used until the last paragraph of
the proof of Proposition \prpref{2mPlus1} cannot be removed. Indeed, it
turns out that no nonassociative Jordan loop of order $5$ or $9$ exists.
In the next section we give an easy proof of this assertion for order $5$.
We were not able to find a short, human argument for order $9$. However,
a finite model builder, such as Mace4 \cite{mace4}, shows that no such loop
exists. Theorem \thmref{Main} will then be  established.

\section{Powers in Jordan loops}
\seclabel{Powers}

Let $Q$ be a loop and $c\in Q$. We use the right-associated convention for
powers: $c^0 = 1$, $c^n = c c^{n-1}$ for $n \geq 1$. For $n > 0$, we say that
$c^n$ is \emph{well-defined} if the value of $c\cdots c$ is independent of how
the $n$ factors are parenthesized. This has an obvious formal characterization,
which could also be taken as an inductive definition.

\begin{lemma}
\lemlabel{well-def}
Let $Q$ be a loop and $c\in Q$. For $n > 0$, $c^n$ is well-defined
if and only if, for each $0 < k < n$, $c^k$ is well-defined and
$c^k c^{n-k} = c^n$.
\end{lemma}

\begin{lemma}\lemlabel{Fifth}
Let $Q$ be a Jordan loop and $c\in Q$. Then $c^n$ is well-defined whenever
$1\le n\le 5$.
\end{lemma}

\begin{proof}
Commutativity alone shows that $c^3$ is well-defined. By \eqnref{Jordan},
$c^2\cdot cc = c^2 c\cdot c$, and hence $c^4$ is well-defined. Finally,
$c^2 c^3 = c^2\cdot c^2 c = c^2 c^2\cdot c$ by \eqnref{Jordan}, and so
$c^5$ is well-defined.
\end{proof}

\begin{lemma}
\lemlabel{Monogenic} Let $Q$ be a loop of order $n$ and let $c\in Q$. If $c^m$
is well-defined for every $1\le m\le n-1$ then $\sblp{c}$ is a cyclic group of
order $k$, and $k=n$ whenever $k>\lfloor n/2\rfloor$.
\end{lemma}

\begin{proof}
The $(n+1)$-tuple $(1,c,\dots,c^{n-1},c\cdot c^{n-1})$ contains a repetition.
Upon canceling $c$ on the left as many times as needed, we conclude that there
is $1\le k\le n$ such that $1$, $c$, $\dots$, $c^{k-1}$ are distinct and
$cc^{k-1}=1$. (We write $cc^{k-1}$ rather than $c^k$ because $k=n$ could occur
and we do not know yet if $c^n$ is well-defined.) Let
$\{a_1,\dots,a_{n-k}\}=Q\setminus\{1,c,\dots,c^{k-1}\}$.

Assume that $k<n$. We have $c^ic^j = c^{(i+j)\mod k}$ for every $0\le i\le
k-1$, $0\le j \le n-k$, as $c^k=1$. If $2k\le n+1$, it immediately follows that
$\langle c\rangle$ is a cyclic group of order $k$. If $2k\ge n$, the nonempty
$(n-k)\times (n-k+1)$ block with rows labeled by $a_1$, $\dots$, $a_{n-k}$ and
columns labeled by $1$, $c$, $\dots$, $c^{n-k}$ must contain only the elements
$a_1$, $\dots$, $a_{n-k}$, a contradiction.

Now assume that $k=n$, the elements of $Q$ are listed as $1$, $c$, $\dots$,
$c^{n-1}$, and $cc^{n-1}=1$. We prove by induction on $i$ that $c^ic^j =
c^{i+j\mod n}$ for every $0\le i$, $j<n$. There is nothing to show when $i=0$,
and the statement also holds for $i=1$ thanks to $cc^{n-1}=1$. Assume that the
statement is true for $i-1\ge 1$, and consider the row labeled by $c^i$, as
visualized below:
\begin{displaymath}
    \begin{array}{c|cccccccc}
    &1&c&\dots&c^{n-i-1}&c^{n-i}&\dots&c^{n-2}&c^{n-1}\\
    \hline
    c^i&c^i&c^{i+1}&\dots&c^{n-1}&?&?&?&?
    \end{array}
\end{displaymath}
We certainly have $c^ic^{n-1}\in\{1,c,\dots,c^{i-1}\}$, and since each of $1$,
$c$, $\dots$, $c^{i-2}$ already occurs in the last column, we in fact have
$c^ic^{n-1} = c^{i-1}$. Proceeding from right to left, we complete the row as
claimed. In particular, $c^ic^{n-i}=1$ since no other power of $c$ is available
at that point.
\end{proof}

\begin{corollary}
\corlabel{assoc5}
Every Jordan loop of order $5$ is a cyclic group.
\end{corollary}

\begin{proof}
Let $Q$ be such a Jordan loop, and fix $c\in Q$ with $c\neq 1$.
By Lemma \lemref{Fifth}, $c^k$ is well-defined for $1\leq k \leq 5$.
By Lemma \lemref{Monogenic}, $\sblp{c}$ is a cyclic group of order
$2$ or $5$. But by Corollary \corref{2Odd}, $c$ cannot have order $2$.
Thus $Q = \sblp{c}$.
\end{proof}

It is not difficult to use Lemma \lemref{Monogenic} to classify Jordan
loops of orders $6$ and $7$. It turns out that there is only one
nonassociative Jordan loop of order $6$, and it is the commutative loop
of exponent $2$ constructed in \S\secref{Even}.
There are only two nonassociative Jordan loops of order $7$, and
each element of those loops has order $3$. One of them is given by
the amalgam construction of \S\secref{Amalgam}, and the other is
given by the construction of \S\secref{Simple}. Details will appear
elsewhere \cite{Pula}.

The case of Jordan loops of order $9$ is trickier. It turns out that
the only two such loops are the two groups. A complete, and rather long,
human proof of this will appear in \cite{Pula}. In the meantime, it is easy to
use a finite model builder such as Mace4 \cite{mace4} to search
exhaustively for Jordan loops of small orders, and such a search
quickly shows that there are no nonassociative models of order $9$.
With this caveat in place, this completes the proof of Theorem \thmref{Main}.

It is an interesting problem to determine which powers of elements
in Jordan loops are well-defined. We give a complete answer for
those powers greater than $5$ which are neither odd primes nor
powers of $2$.

\begin{theorem}
\thmlabel{powers}
Let $m,n$ be integers with $m\geq 2$, $n \geq 3$, and $n$ odd. Then there exists
a $1$-generated Jordan loop $Q$ with a generator $c$ such that $c^k$ is well-defined for
$0\leq k < mn$, but $c^{pn} c^{(m-p)n} \neq c^{mn}$ for $0 < p < m$.
\end{theorem}

\begin{proof}
Let $s\ge m+2$ be the smallest positive integer relatively prime to $n$.

For an integer $\ell$, we denote by $\bar{\ell}$ its image in
$\mathbb{Z}_n = \mathbb{Z}/n\mathbb{Z}$ and by $[\ell]$ the
image in $\mathbb{Z}_s$.
Then $[n]$ is a generator of
$\mathbb{Z}_s = \{[0], [n], \ldots , (s-1)[n] \}$.

We define a permutation $\phi$ on $\mathbb{Z}_s$ by
\[
\phi(i [n]) = \begin{cases}
i [n], & \text{if } 0 \leq i \leq m - 1 \\
(s + m - i - 1)[n], & \text{if } m \leq i \leq s - 1
\end{cases}
\]
Using $\phi$, we define a new operation $\circ$ on $\mathbb{Z}_s$ by
\[
[k] \circ [\ell] = \phi^{-1} (\phi([k]) + \phi([\ell]))
\]
for all $k,l\in \mathbb{Z}$. (The assumption that $s \geq m + 2$
is used here to guarantee that $\phi(m[n]) \neq m[n]$.)
Then $(\mathbb{Z}_s,\circ)$ is an isomorphic copy of
$(\mathbb{Z}_s,+)$.

Now let $Q = \mathbb{Z}_n \times \mathbb{Z}_s$, and define
a binary operation $\cdot$ on $Q$ as follows. For $a,b,u,v\in \mathbb{Z}$,
set
\[
(\bar{a},[u])\cdot (\bar{b},[v]) =
\begin{cases}
(\overline{a + b}, [u + v]) & \text{if } a\not\equiv 0 \text{ or } b\not\equiv 0 \mod n \\
(\bar{0}, [u]\circ [v]) & \text{if } a\equiv b \equiv 0 \mod n
\end{cases}
\]

It is clear that $(Q,\cdot)$ is a commutative loop with neutral element $(0,0)$.
We claim that $(Q,\cdot)$ is, in fact, a Jordan loop, that is, it satisfies
\[
(\bar{a},[u])^2 \cdot (\bar{b},[v])(\bar{a},[u])
= (\bar{a},[u])^2 (\bar{b},[v])\cdot (\bar{a},[u])
\]
for all $a,b,u,v\in \mathbb{Z}$.

Firstly, assume $a\not\equiv 0 \mod n$. Then $(\bar{a},[u])^2 =
(\bar{2a},[2u])$. Since $n$ is odd, $2a \not\equiv 0 \mod n$, and so
\begin{align*}
(\bar{a},[u])^2 (\bar{b},[v])\cdot (\bar{a},[u])
&= (\overline{2a+b},[2u+v])(\bar{a},[u]) = (\overline{3a+b},[3u+v]) \\
&= (\bar{2a},[2u])\cdot (\bar{b},[v])(\bar{a},[u])
= (\bar{a},[u])^2\cdot (\bar{b},[v])(\bar{a},[u])\,.
\end{align*}

Next, assume $a \equiv 0\mod n$. Then $(\bar{0},[u])^2 =
(\bar{0},[u]\circ [u])$. If $b\not\equiv 0\mod n$, then
\begin{align*}
(\bar{0},[u])^2 \cdot (\bar{b},[v])(\bar{0},[u])
&= (\bar{0},[u]\circ [u])\cdot (\bar{b},[v+u])
= (\bar{b}, [([u]\circ [u]) + v + u]) \\
&= (\bar{b},[([u]\circ [u])+v])(\bar{0},[u])
= (\bar{0},[u])^2 (\bar{b},[v]) \cdot (\bar{0},[u])\,.
\end{align*}
The identity also holds if $b \equiv 0\mod n$, for
then all calculations in the second component occur in the cyclic group
$(\mathbb{Z}_s,\circ)$. This completes the verification that
$(Q,\cdot)$ is a Jordan loop.

Set $c = (\bar{1}, [1]) \in Q$. Then for all $k > 0$,
$c^k = (\bar{k}, [k])$. (Recall that $c^k$ is the
right associated product $c(c\cdots c)$.)
We claim that $c$ generates $Q$. Indeed, for $a,b\in \mathbb{Z}$,
choose $p,q,r,t\in \mathbb{Z}$ such that
$p n + q s = a$ and $r n + t s = b$. Then
$c^{q s} c^{r n} = ( \overline{q s}, [0] )( \bar{0}, [r n] )
= ( \overline{q s}, [r n] ) = ( \bar{a} , [b] )$.

Next we wish to show  that $c^k$ is well-defined for $0\leq k < mn$. Thus we
must show that $c^i c^j = c^{i+j}$ for $0\leq i + j \leq mn-1$.

If at least one of $i,j$ is not divisible by $n$, then
$c^i c^j = (\bar{i},[i])(\bar{j},[j]) =
(\overline{i+j},[i+j]) = c^{i+j}$. For the other case,
assume that $i = pn$ and $j = qn$, where
$0\leq p,q,p + q \leq m-1$. Then
\[
[i] \circ [j] = \phi^{-1}(\phi(p[n]) + \phi(q[n]))
= \phi^{-1}((p+q)[n]) = (p+q)[n] = [i+j]\,.
\]
Therefore $c^i c^j = (\bar{0},[i])(\bar{0},[j])
= (\bar{0},[i]\circ [j]) = (\bar{0},[i] + [j])
= c^{i + j}$ in this case as well.

Finally, suppose $0 < p < m$. Then $m-p \leq m-1$, and so
$\phi^{-1}(\phi(p[n]) + \phi((m-p)[n])) =
\phi^{-1}(p[n] + (m-p)[n]) = \phi^{-1}(m[n]) = (s-1)[n]$.
Therefore,
\begin{align*}
c^{pn} c^{(m-p)n} &= (\bar{0},p[n])(\bar{0},(m-p)[n]) = (\bar{0},p[n]\circ (m-p)[n]) \\
&= (\bar{0},(s-1)[n]) \neq (\bar{0},m[n]) = c^{mn}\,,
\end{align*}
as claimed.
\end{proof}

To illustrate the construction of Theorem \thmref{powers},
consider the particular case $n = 3$, $m = 2$.
In this case, $s = 4$, and so the order of the loop will be $12$.
The permutation $\phi$ is given by $\phi([0]) = [0]$,
$\phi([3]) = [3]$, $\phi([2]) = \phi(2[3]) = 3[3] = [1]$,
$\phi([1]) = \phi(3[3]) = 2[3] = [2]$. The operation $\circ$
on $\mathbb{Z}_4$ is given by
\[
\begin{array}{c|cccc}
\circ & \lbrack 0\rbrack& \lbrack 1\rbrack& \lbrack 2\rbrack& \lbrack 3\rbrack\\
    \hline
\lbrack 0\rbrack & \lbrack 0\rbrack& \lbrack 1\rbrack& \lbrack 2\rbrack& \lbrack 3\rbrack \\
\lbrack 1\rbrack & \lbrack 1\rbrack& \lbrack 0\rbrack& \lbrack 3\rbrack& \lbrack 2\rbrack \\
\lbrack 2\rbrack & \lbrack 2\rbrack& \lbrack 3\rbrack& \lbrack 1\rbrack& \lbrack 0\rbrack \\
\lbrack 3\rbrack & \lbrack 3\rbrack& \lbrack 2\rbrack& \lbrack 0\rbrack& \lbrack 1\rbrack
\end{array}
\]
We represent the loop element $(\bar{a},[u])$ by the nonnegative integer
$a + 3u$. Then the loop multiplication table is as follows.
\[
\begin{array}{ccc|ccc|ccc|ccc}
0& 1& 2&            3& 4& 5&            6& 7& 8&             9&10&11 \\
1& 2& 0&            4& 5& 3&            7& 8& 6&            10&11& 9 \\
2& 0& 1&            5& 3& 4&            8& 6& 7&            11& 9&10 \\
\hline
3& 4& 5&            0& 7& 8&            9&10&11&            6& 1& 2 \\
4& 5& 3&            7& 8& 6&            10&11& 9&           1& 2& 0 \\
5& 3& 4&            8& 6& 7&            11& 9&10&           2& 0& 1 \\
\hline
6& 7& 8&            9&10&11&            3& 1& 2&            0& 4& 5\\
7& 8& 6&        10&11& 9&           1& 2& 3&            4& 5& 6\\
8& 6& 7&        11& 9& 10&          2& 3& 4&            5& 6& 7\\
\hline
 9&10&11&           6& 1& 2&            0& 4& 5&        3& 7& 8\\
10&11& 9&         1& 2& 0&          4& 5& 6&        7& 8& 6\\
11& 9& 10&      2& 0& 1&            5& 6& 7&        8& 6& 7
\end{array}
\]
The generator $c$ of Theorem \thmref{powers} is $4$, and
we see that $4(4\cdot 4)\cdot 4(4\cdot 4) = 3 \neq 6 =
4\cdot 4(4\cdot 4(4\cdot 4))$.

Whether $x^m$ where $m > 6$ is an odd prime or a power of $2$ is well-defined
is at present unclear. For instance, in a Jordan loop, if $c$ is an element
such that $c^6$ is well-defined, then both $c^7$ and $c^8$ are well-defined. If
a Jordan loop $Q$ satisfies $x^m$ is well-defined for \emph{all} $x\in Q$ for
$0 \leq m \leq 10$, then $x^{11}$ is also well-defined for all $x\in Q$. We do
not know of an example of a Jordan loop in which there exists an element $c$
with $c^m$ well-defined for $0\leq m\leq 10$, but such that $c^{11}$ is not
well-defined.

\begin{problem}
Let $Q$ be a Jordan loop, and
let $m > 6$ be an integer which is either an odd prime
or a power of $2$.
\begin{enumerate}
\item For a fixed $c\in Q$, if $c^k$ is well-defined for
$0\leq k < m$, must $c^m$ be well-defined?
\item If $x^k$ is well-defined for all $x\in Q$ and for
$0\leq k < m$, must $x^m$ be well-defined?
\end{enumerate}
\end{problem}

\begin{problem}
For which $n$ does there exist a $1$-generated Jordan loop
of order $n$ which is not power associative? In particular,
is there a $1$-generated Jordan loop of prime order or of
order a power of $2$ which is not power associative?
\end{problem}

\section{An infinite family of nonassociative simple Jordan loops}
\seclabel{Simple}

In this section we construct an infinite family of nonassociative simple Jordan
loops. Although the construction could be viewed as another generalized amalgam,
we opt for a more direct description.

Let $A=(A,\cdot)$ be a commutative loop of order $2^n-1\ge 1$, and let
\begin{displaymath}
    B=B_n=\{[x_1,\dots,x_n]=[x_i];\;x_i\in\{0,1\}\}
\end{displaymath}
be the binary hypercube of size $2^n$. Label elements of $A$ by $(x_i)$, where
$[x_i]$ is any element of $B$ different from $[1,\dots,1]$, and let
$(0,\dots,0)$ be the neutral element of $A$. (As a mnemonic device, loop
elements are enclosed in round parentheses, hypercube elements in brackets.)

Define an operation $\circ$ on the disjoint union $A\cup B$ by
\begin{align*}
    (x_i)\circ (y_i) &= (x_i)\cdot (y_i),\\
    (x_i)\circ [y_i] &= [x_i\oplus y_i],\\
    [x_i]\circ (y_i) &= [x_i\oplus y_i],\\
    [x_i]\circ [y_i] &= \left\{
        \begin{array}{ll}
            [x_i'],&\text{ if $[x_i]=[y_i]$},\\
            ((x_i\oplus y_i)'),&\text{ otherwise},
        \end{array}\right.
\end{align*}
where $\oplus$ is the usual binary addition, and $'$ is the binary complement.

Note that $[(x_i\oplus y_i)']$ is different from $[1,\dots,1]$ when $[x_i]\ne
[y_i]$, and hence $\circ$ is indeed an operation on $A\cup B$, clearly
commutative. We will denote the resulting groupoid by $\jordan{A}$.

Before we analyze the construction, it is once again worth visualizing it in
terms of a multiplication table.

Here is $\mathcal J(1)$ for the trivial one-element loop $1$:
\begin{displaymath}
    \begin{array}{c||c|cc}
        \circ   &   (0) &   [0] &   [1]\\
        \hline\hline
        (0)     &   (0) &   [0] &   [1]\\
        \hline
        {[0]}     &   [0] &   [1] &   (0)\\
        {[1]}     &   [1] &   (0) &   [0]
    \end{array}.
\end{displaymath}
Of course, this is just the cyclic group of order $3$, hence again a
commutative loop. Upon relabeling the elements of $\jordan{1}$ as $(0,0)$,
$(0,1)$, $(1,0)$, we obtain $\jordan{\jordan{1}}$:
\begin{displaymath}
    \begin{array}{c||ccc|cccc}
    \circ& (0,0)& (0,1)& (1,0)& [0,0]& [0,1]& [1,0]& [1,1]\\
    \hline\hline
    (0,0)& (0,0)& (0,1)& (1,0)& [0,0]& [0,1]& [1,0]& [1,1]\\
    (0,1)& (0,1)& (1,0)& (0,0)& [0,1]& [0,0]& [1,1]& [1,0]\\
    (1,0)& (1,0)& (0,0)& (0,1)& [1,0]& [1,1]& [0,0]& [0,1]\\
    \hline
    {[0,0]}& [0,0]& [0,1]& [1,0]& [1,1]& (1,0)& (0,1)& (0,0)\\
    {[0,1]}& [0,1]& [0,0]& [1,1]& (1,0)& [1,0]& (0,0)& (0,1)\\
    {[1,0]}& [1,0]& [1,1]& [0,0]& (0,1)& (0,0)& [0,1]& (1,0)\\
    {[1,1]}& [1,1]& [1,0]& [0,1]& (0,0)& (0,1)& (1,0)& [0,0]
    \end{array}
\end{displaymath}
We can therefore think of the multiplication table for $\jordan{A}$ as
consisting of four quadrants:
\begin{enumerate}
\item[-] the top left $2^n-1$ times $2^n-1$ quadrant is the multiplication
table of $A$,

\item[-] the top right $2^n-1$ times $2^n$ quadrant is the canonical
multiplication table of the elementary abelian group of order $2^n$ with the
last row removed,

\item[-] the bottom left $2^n$ times $2^n-1$ quadrant is the canonical
multiplication table of the elementary abelian group of order $2^n$ with the
last column removed,

\item[-] the bottom right $2^n$ times $2^n$ quadrant is the canonical
multiplication table of the elementary abelian group of order $2^n$, except
that the elements are listed in reversed order (this is because of the
complement in the definition), and that the elements of $B_n$ are superimposed
on the main diagonal, again in reversed order.
\end{enumerate}

A short reflection shows that the multiplication table of $\jordan{A}$ is a
Latin square. We have obtained:

\begin{proposition} Let $A$ be a commutative loop of order $2^n-1$, and let $\jordan{A}$ be
defined as above. Then $\jordan{A}$ is a commutative loop of order $2^{n+1}-1$
with subloop $A$.
\end{proposition}

\begin{proposition}\prplabel{IsJordan}
Let $A$ be a Jordan loop. Then $\jordan{A}$ is a Jordan
loop.
\end{proposition}
\begin{proof}
We must verify $(x\circ x)\circ (y\circ x) = ((x\circ x)\circ y)\circ x$ for
every $x$, $y\in\jordan{A} = A\cup B$. When $x\in A$, $y\in A$, there is
nothing to show since $A$ is Jordan and $\circ$ coincides with the
multiplication in $A$ on $A\times A$. For the rest of the proof, let $x=(x_i)$
or $x=[x_i]$, and $y=(y_i)$ or $y=[y_i]$.

Assume that $x\in A$ and $y\in B$. Let $(z_i) = x\circ x$. Then $(x\circ
x)\circ(y\circ x) = ((x_i)\circ(x_i))\circ ([y_i]\circ(x_i)) = (z_i)\circ
[y_i\oplus x_i] = [z_i\oplus y_i\oplus x_i]$. On the other hand, $((x\circ
x)\circ y)\circ x = (((x_i)\circ(x_i))\circ[y_i])\circ(x_i) =
((z_i)\circ[y_i])\circ(x_i) = [z_i\oplus y_i]\circ (x_i) = [z_i\oplus y_i\oplus
x_i]$.

Assume that $x\in B$ and $y\in A$. Then $(x\circ x)\circ(y\circ x) =[x_i']\circ
[y_i\oplus x_i]$. Note that $x_i' = y_i\oplus x_i$ for every $i$ if and only if
$y_i=1$ for every $i$. But this cannot happen since $y\in A$. Thus $[x_i']\circ
[y_i\oplus x_i] = (x_i'\oplus y_i\oplus x_i) = (y_i)$. On the other hand,
$((x\circ x)\circ y)\circ x = ([x_i']\circ (y_i))\circ [x_i] = [x_i'\oplus
y_i]\circ [x_i]$. For the same reason as above we cannot have $x_i'\oplus
y_i=x_i$ for every $i$, and thus $[x_i'\oplus y_i]\circ[x_i] = (y_i)$.

Finally assume that $x\in B$ and $y\in B$. Then $(x\circ x)\circ(y\circ x)$ is
equal to
\begin{displaymath}
    [x_i']\circ ([y_i]\circ[x_i]) = \left\{
        \begin{array}{ll}
            [x_i']\circ(y_i\oplus x_i) = [y_i],&\text{ if $[x_i]\ne[y_i]$},\\
            {[x_i']}\circ[x_i'] = [x_i'']=[x_i]=[y_i],&\text{ otherwise}.
        \end{array}\right.
\end{displaymath}
On the other hand, $((x\circ x)\circ y)\circ x$ is equal to
\begin{displaymath}
    ([x_i']\circ [y_i])\circ[x_i] = \left\{
        \begin{array}{ll}
            (x_i'\oplus y_i)\circ[x_i] = [y_i],&\text{ if $[x_i']\ne[y_i]$},\\
            {[x_i'']}\circ[x_i] = [x_i]\circ[x_i]=[x_i']=[y_i],&\text{ otherwise}.
        \end{array}\right.
\end{displaymath}
This finishes the proof.
\end{proof}

\begin{lemma}\lemlabel{Order3}
Let $A$ be a commutative loop. Then every element of $\jordan{A}\setminus A$
has order $3$. In particular, when $A$ is power associative then so is
$\jordan{A}$.
\end{lemma}
\begin{proof}
We must show that $(x\circ x)\circ x = x\circ (x\circ x) = (0)$ for every
$x=[x_i]\in\jordan{A}\setminus A$. We have $(x\circ x)\circ x =
[x_i']\circ[x_i] = ((x_i'\oplus x_i)') = (1') = (0)$, and similarly for $x\circ
(x\circ x)$.
\end{proof}

\begin{lemma}\lemlabel{NotLeftAlt}
Let $A$ be a commutative loop of order $2^n-1>1$. Then $\jordan{A}$ is not left
alternative, and hence not diassociative.
\end{lemma}
\begin{proof}
The commutative loop $A$ cannot be of exponent $2$ by Corollary \corref{2Odd}.
Let $x=(x_i)\in A$ be an element satisfying $x^2\ne 1$. Let $y=[0]\in B$. Then
$x\circ(x\circ y) = (x_i)\circ [x_i] = [0]$. Now, $(z_i)= x\circ x$ is a
nonidentity element of $A$, and thus $(x\circ x)\circ y = [z_i]\ne[0]$.
\end{proof}

If $Q$ is a loop and $x\in Q$, denote by $L_x$, $R_x$ the left and right
translations by $x$ in $Q$. Recall the \emph{left inner mappings} $L(x,y) =
L^{-1}_{yx}L_yL_x$, the \emph{right inner mappings} $R(x,y) =
R_{xy}^{-1}R_yR_x$, and the \emph{middle inner mappings}, or
\emph{conjugations} $T(x) = R_x^{-1}L_x$.

A subloop $S$ of a loop $Q$ is \emph{normal}, $S\unlhd Q$, if $S$ is invariant
under all inner mappings of $Q$. A loop $Q$ is \emph{simple} if it has no
normal subloops except for $\{1\}$ and $Q$.

\begin{proposition}\prplabel{Simple}
Let $A$ be a simple commutative loop.
Then $\jordan{A}$ is simple as well.
\end{proposition}
\begin{proof}
The statement holds when $|A|=1$, so suppose that $|A|=2^n-1>1$. Let $S$ be a
nontrivial normal subloop of $\jordan{A}=A\cup B$. Then $S\cap A$ is a normal
subloop of $A$. Since $A$ is simple, we have $S\cap A=1$ or $S\cap A=A$. The
latter case implies that $|S|\ge 2^n-1$. Since $|S|$ has to be a divisor of
$|\jordan{A}|=2^{n+1}-1$, it follows that $S=\jordan{A}$, a contradiction.

So suppose that $S\cap A=1$, and let $x\in S\setminus A$. By Lemma
\lemref{Order3}, $|x|=3$.

If $S$ contains $\langle x\rangle$ properly, let $y\in S\setminus(A\cup\langle
x\rangle)$. Then $x\circ y = [x_i]\circ [y_i] = ((x_i\oplus y_i)')$. Note that
$[y_i]=y\ne x^2=[x_i']$, and thus $x\circ y = ((x_i\oplus y_i)')\ne (0)$,
contradicting $S\cap A=1$.

We can therefore assume that $S=\langle x\rangle$ for some $x\in S\setminus A$.
Recall that $S$ is closed under all left inner mappings of $\jordan{A}$. Thus
for every $u$, $v\in\jordan{A}$ there must be $z\in S$ such that $u\circ(v\circ
x) = (u\circ v)\circ z$. We show that this is not the case.

First assume that $x\ne[0]$ and $x\ne[1]$. Then $[0]\circ([0]\circ[x_i]) =
[0]\circ(x_i') = [x_i']$. But $[0]\circ[0]=[1]\ne[x_i']$,
$([0]\circ[0])\circ[x_i] = [1]\circ[x_i] = (x_i)\ne[x_i']$, and
$([0]\circ[0])\circ([x_i]\circ[x_i]) = [1]\circ[x_i'] = (x_i')\ne[x_i']$. Now
assume that $x=[0]$, so $x^2=[1]$. Pick $y\in A$ such that $y^2\ne 1$. Then
$(y_i)\circ((y_i)\circ[0]) = (y_i)\circ[y_i] = [0]$, but
$(y_i)\circ(y_i)=(z_i)\ne[0]$, $((y_i)\circ(y_i))\circ[0] =(z_i)\circ[0] =
[z_i]\ne [0]$, and $((y_i)\circ(y_i))\circ[1] = (z_i)\circ[1] = [z_i']$. Note
that $[z_i']\ne [0]$, else $[z_i]=[1]$, $(z_i)=(1)$, a contradiction with
$(z_i)\in A$. The remaining case $x=[1]$ is analogous.
\end{proof}

An infinite family of nonassociative simple Jordan loops is now at hand. Let
$A_0$ be the trivial group, and $A_{i+1} = \jordan{A_i}$ for every $i\ge 0$. By
Proposition \prpref{IsJordan}, every $A_i$ is a Jordan loop. By Proposition
\prpref{Simple}, every $A_i$ is simple. By Lemma \lemref{NotLeftAlt}, $A_i$ is
nonassociative for every $i\ge 2$.

Our results also imply that $\jordan{\mathbb{Z}_p}$ is a nonassociative simple
Jordan loop whenever $p>2$ is a Mersenne prime.

\section{Acknowledgment}

We thank the referees for several helpful comments, including the
succinct algebraic description given in the proof of Lemma \lemref{IdempQuasi}.

\bibliographystyle{plain}

\begin{thebibliography}{99}

\bibitem{Bruck44} R.~H.~Bruck,
Some results in the theory of linear non-associative algebras,
\textit{Trans. Amer. Math. Soc.} \textbf{56} (1944), 141-199.

\bibitem{Bruck}
R.~H.~Bruck,
\textit{A Survey of Binary Systems},
\textit{Ergebnisse der Mathematik und Ihrer Grenzgebiete}, New Series, Volume
\textbf{20}, Springer, 1958.

\bibitem{CPS} O.~Chein, H.~O.~Pflugfelder and J.~D.~H.~Smith (eds.),
\textit{Quasigroups and Loops: Theory and Applications}, Sigma Series in Pure
Mathematics \textbf{8}, Heldermann, 1990.

\bibitem{ColbournRosa}
C.~J.~Colbourn and A.~Rosa,
\textit{Triple systems}, Oxford Mathematical Monographs,
Oxford Science Publications, Clarendon Press, Oxford 1999.

\bibitem{Foguel} T.~Foguel,
Amalgam of a loop over an idempotent quasigroup,
\textit{Quasigroups and Related Systems} \textbf{13} (2005), 99--104.

\bibitem{GoodaireKeeping}
E.~G.~Goodaire and R.~G.~Keeping,
Jordan loops and loop rings, preprint.

\bibitem{Kirkman}
T.~P.~Kirkman,
On a problem in combinations,
\textit{Cambridge and Dublin Math. J.} \textbf{2} (1897), 191--204.

\bibitem{McCrimmon}
K.~McCrimmon, \textit{A Taste of Jordan Algebras}, Universitext, Springer, 2004.

\bibitem{mace4}
W.~W.~McCune,
\textit{Mace4 Reference Manual and Guide},
Tech. Memo ANL/MCS-TM-264, Mathematics and Computer Science Division,
Argonne National Laboratory, Argonne, IL, August 2003.
\texttt{http://www.cs.unm.edu/\symbol{126}mccune/mace4/}

\bibitem{pflug} H.~O.~Pflugfelder,
\textit{Quasigroups and Loops: Introduction}, Sigma Series in Pure
Mathematics \textbf{7}, Heldermann Verlag, 1990.

\bibitem{Pula} K.~Pula,
\textit{Powers of elements in Jordan loops},
Commentationes Mathematicae Universitatis Carolinae, to appear.

\end{thebibliography}

\end{document}